\documentclass{amsart}

\usepackage{amssymb,amsmath,amsfonts,latexsym,amsthm,graphicx} 

\usepackage[ngerman,english]{babel}

\usepackage[latin1]{inputenc}
\usepackage{hyperref}

\usepackage{color}

\newtheorem{theorem}{Theorem}[section]

\newtheorem{Corollary}[theorem]{Corollary}

\newtheorem{Lemma}[theorem]{Lemma}

\theoremstyle{definition}

\newtheorem{Remark}[theorem]{Remark}

\DeclareMathOperator*{\esssup}{ess\,sup}

\usepackage{tikz}

\usepackage{tikz-cd}

\newcommand{\C}{\mathbb{C}} 
\newcommand{\R}{\mathbb{R}} 
\newcommand{\Z}{\mathbb{Z}} 
\newcommand{\N}{\mathbb{N}} 

\newcommand{\wh}{\widehat}
\newcommand{\vp}{\varphi}
\newcommand{\eps}{\varepsilon}

\title{A Fourier integral formula for logarithmic energy}
\author{L.\ Frerick}
\address{University of Trier, FB IV, Mathematics  D-54286 Trier, Germany}
\email{frerick@uni-trier.de}

\author{J.\ M\"uller} 
\address{University of Trier, FB IV, Mathematics  D-54286 Trier, Germany}
\email{jmueller@uni-trier.de}

\author{T.\ Thomaser}
\address{University of Trier, FB IV, Mathematics  D-54286 Trier, Germany}
\email{thomaser@uni-trier.de}

\begin{document}

\begin{abstract}
 A formula which expresses logarithmic energy of Borel measures on $\R^n$  in terms of the Fourier transforms of the measures is established and some applications are given. In addition, using similar techniques a (known) formula for Riesz energy is reinvented. 
\end{abstract}



\maketitle

\section{Introduction and main result}

The notions of logarithmic potential and logarithmic energy play a central role  in potential theory in particular in dimension two as well as in free probability. Recommended introductions are  \cite{Landkof} and \cite{Ransford} for the classical theory and \cite{MS} for the more recent theory of free probability.  
We consider arbitrary dimension $n  \in \N$ and write $\mathcal{M}(\R^n)$ for the set of all complex Borel measures on $\R^{n}$ and $\mathcal{M}_+(\R^n)$ for the subset of nonnegative $\mu$. With $|\mu|$ denoting the total variation of $\mu \in \mathcal{M}(\R^n)$ and $|x|$ denoting the euclidean norm of $x \in \R^n$
we write  
$$
p_\mu(x):=\displaystyle\int \ln\left(\frac{1}{\vert x-y\vert}\right) \, d\mu(y) \in \C \cup\{\pm \infty\}
$$
if $y \mapsto \ln|x-y|$ is integrable with respect to $|\mu|$ or $\mu \in \mathcal{M}_+(\R^n)$ and the integral exists in $[-\infty, \infty]$. 
 Moreover, we write $\mathcal{M}_{0}(\R^{n})$ for set of all $\mu \in \mathcal{M}(\R^n)$ with 
$$
\displaystyle\int\displaystyle\int \left\vert\ln\left(\vert x-y\vert\right)\right\vert \, d\vert\mu\vert(y) \, d\vert\mu\vert(x)<+\infty
$$
and $\mathcal{M}_{0,+}(\R^n)$ for the set of all $\mu \in \mathcal{M}_+(\R^n)$  with the property that the integral
$
\int\int \ln\vert x-y\vert \, d\mu(y) \, d\mu(x)
$
exists in  $[-\infty,+\infty]$. 
For  $\mu\in\mathcal{M}_{0}(\R^{n}) \cup \mathcal{M}_{0,+}(\R^n)$ the logarithmic potential $p_\mu$ of $\mu$  exists
$\mu$-almost everywhere  and 
$$
I(\mu):= \displaystyle\int\displaystyle\int \ln\left(\frac{1}{\vert x-y\vert}\right) \, d\mu(y) \, d\bar{\mu}(x)= \displaystyle\int p_{\mu}(x) \, d\bar{\mu}(x) 
$$
is called the logarithmic energy of $\mu$. By definition, $I(\mu)$ is a complex number for $\mu \in  \mathcal{M}_{0}(\R^n)$ and $I(\mu)\in [-\infty, \infty]$ for $ \mu \in \mathcal{M}_{0,+}(\R^n)$.
If $\mu \ge 0$ with 
$$\int\ln(1+\vert x\vert) \, d\mu(x)<+\infty$$
 then  $I(\mu)>-\infty$ thanks to $
\vert x-y\vert\leq (1+\vert x\vert)(1+\vert y\vert)$. In particular, every $\mu\ge 0$ with compact support lies in $\mathcal{M}_{0, +}(\R^{n})$.

Let
\[
\mathbb{S}^{n-1}:=\{ x \in \R^n : |x|=1\}
\]
denote the $(n-1)$-sphere, $\sigma_{n-1}$ the surface measure of $\mathbb{S}^{n-1}$  and $\omega_{n-1}$ the corresponding area, given by $2\pi^{n/2}/\Gamma(n/2)$. In particular,  $\mathbb{S}:=\mathbb{S}^1$ is the unit circle in $\C=\R^2$. Writing $a_k:=\int_\mathbb{S} \zeta^k \, d \mu(\zeta)$ for Fourier coefficients of a complex measure $\mu$ supported on $\mathbb{S}$, a basic fact is that the logarithmic energy can be expressed  by 
\begin{equation}\label{discrete}
2I(\mu)=\sum\limits_{k \in \Z \setminus \{0\}} \dfrac{\vert a_k\vert^{2}}{|k|}\, .
\end{equation}
 This is an identity if $\mu \in \mathcal{M}_{0}(\R^2)$ or $\mu \ge 0$ and serves as definition of $I(\mu)$ for arbitrary complex $\mu$ (cf.\ \cite{Ka}, p.\ 35, \cite{HR2}).  In particular, in all cases $I(\mu) \in [0, \infty]$. The identity can be seen by expanding
\[
 2\ln\left(\frac{1}{\vert z-rw\vert}\right) = 2\, {\rm Re}\log\left(\frac{1}{1-rz\overline{w}}\right) = \sum_{k \in \Z \setminus \{0\}} \frac{r^{|k|}}{|k|}z^k\overline{w}^k
\] 
for $0<r<1$ and taking the limit $r \to 1$ (cf.\ \cite{El-Fallah}, Section 2.4). 
 It might come as a surprise that $a_0=\mu(\mathbb{S})$ does not enter the scene. This can be explained by the fact that the logarithmic potential of the arclength measure $\sigma_1$  vanishes on $\mathbb{S}$ and, as a consequence,  $I(\mu)$ equals $I(\mu + \lambda \sigma_1)$ for arbitrary  $\mu$ and scalars $\lambda$ (cf.\ \cite{Meckes}, p. 119). 

 Our aim is to establish a continuous $n$-dimensional version of a formula for the  logarithmic energy of appropriate measures $\mu$ in terms of the Fourier transform
\[
\widehat{\mu}(\xi):=\int_{\R^n} e^{i\xi \cdot x}\, d\mu(x) \qquad (\xi \in \R^n)
\]
and the generalized hypergeometric function $K_n$, defined by 
\[
K_{n}(z)
=\sum\limits_{k=0}^{\infty} \frac{(-1)^{k}}{(n/2,k)k!} \left(\frac{z}{2}\right)^{2k}\quad (z \in \C)
\]
with the Pochhammer symbol $(\alpha,k):=\Gamma(\alpha +k)/\Gamma(\alpha)$. 
In particular, one obtains
$K_1=\cos$ 
 and in the case  $n=2$, in which the logarithmic potential coincides with the Newtonian, 
$K_2=J_0$, where $J_\alpha$ denotes the Bessel function of first kind and order $\alpha$. 
More generally,  $K_n$ can be expressed in terms of $J_{n/2-1}$ as 
\[
K_n(t)=2^{n/2-1}\Gamma(n/2) J_{n/2-1}(t)/t^{n/2-1} \quad (t >0).
\]  
With these notations, the main result on mutual logarithmic energy reads as follows:  
\begin{theorem} \label{Energy formula several dimensions} 
Let $\mu,\nu\in\mathcal{M}(\R^{n})$ with  
$$
\int\int \left\vert\ln\left(\vert x-y\vert\right)\right\vert \, 
d\vert\mu\vert(y) \, d\vert\nu\vert(x)< \infty
$$
or $\mu, \nu \in\mathcal{M}_{+}(\R^{n})$ with the property that
$\int p_\mu \, d\nu$
exists in  $[-\infty,+\infty]$. Then
\begin{equation}\label{f1}
\omega_{n-1} \int p_\mu\, d\bar{\nu}= \lim\limits_{\substack{\varepsilon\rightarrow 0\\ N\rightarrow\infty}} \ \displaystyle\int\limits_{\varepsilon\leq \vert \xi\vert\leq N} \left((\widehat{\mu}\, \overline{\widehat{\nu}})(\xi)-\mu(\R^{n})\bar{\nu}(\R^{n})K_{n}(\vert \xi\vert)\right) \, \dfrac{d\xi}{\vert \xi\vert^{n}}.
\end{equation}
In particular, for $\mu\in\mathcal{M}_{0}(\R^{n}) \cup \mathcal{M}_{0,+}(\R^{n})$ we have
\begin{equation}\label{f2}
\omega_{n-1}I(\mu)=\lim\limits_{\substack{\varepsilon\rightarrow 0\\ N\rightarrow\infty}} \ \displaystyle\int\limits_{\varepsilon\leq \vert \xi\vert\le N} \left(\vert\widehat{\mu}(\xi)\vert^{2}-\vert\mu(\R^{n})\vert^{2}K_{n}(\vert \xi\vert)\right) \, \dfrac{d\xi}{\vert \xi\vert^{n}}\, .
\end{equation}
\end{theorem}

If $\wh{\mu}\, \overline{\wh{\nu}}-\mu(\R^{n})\bar{\nu}(\R^{n})$ is locally integrable at the origin with respect to $d\xi/|\xi|^n$, then the integrand function on the right hand side in \eqref{f1} has the same property and the double-sided limit reduces to a one-sided $\lim_{N \to \infty}$. This holds if $\wh{\mu}\, \overline{\wh{\nu}}$ is Dini continuous at $0$, which is  in particular  the case if $\mu$ and $\nu$ have  compact support.


We will frequently use the following substitution rule for the $n$-dimensional Lebesgue measure, a  proof of which can be found e.g.\ in \cite{Folland}, p.\ 78: If $f:\R^n \to\C$ is Lebesgue integrable, then
\begin{equation} \label{Lemma sphere}
\displaystyle\int_{\R^{n}} f(\xi) \, d\xi= \displaystyle\int_{0}^{\infty}\displaystyle\int_{\mathbb{S}^{n-1}} f(r\zeta) \, d\sigma_{n-1}(\zeta) \, r^{n-1} \, dr.
\end{equation} 
Due to the fact that $J_{\alpha}(t)=O(t^{-1/2})$ as $0< t\rightarrow \infty$ for $\alpha\geq 0$ (see \cite{Olver}, 10.17.3), the function $t \mapsto K_n(t)/t$ is absolutely integrable  at $+\infty$ for $n \ge 2$ and improperly integrable  for $n=1$. Hence,
\[
\lim\limits_{N\rightarrow\infty} \displaystyle\int_{1\leq \vert \xi\vert\leq N} \dfrac{K_{n}(\vert \xi\vert)}{\vert \xi\vert^{n}} \, d\xi=\omega_{n-1} \lim\limits_{N\rightarrow\infty} \displaystyle\int_{1}^{N} \dfrac{K_{n}(r)}{r} \, dr 
\]
exists and so,  by monotonicity, for arbitrary complex Borel measures  the limit $N \to \infty$ in \eqref{f2} also exists as value in $(-\infty, \infty]$. This implies that the function $|\wh{\mu}|^2$ is integrable at $\infty$ with respect to $d\xi/|\xi|^n$ if $I(\mu)$ is finite and that
\[
I(\mu):=\frac{1}{\omega_{n-1}} \displaystyle\int \left(\vert\widehat{\mu}(\xi)\vert^{2}-\vert\mu(\R^{n})\vert^{2}K_{n}(\vert \xi\vert)\right) \, \dfrac{d\xi}{\vert \xi\vert^{n}}
\]
extends the definition of logarithmic energy to arbitrary complex Borel measures having compact support. 

Choosing $d\mu(x)=\vp(x)\, dx$, where $\vp$ belongs to the Schwartz space $\mathcal{S}$, we obtain from \eqref{f1} with the Dirac measure $\delta$ at $0$ and $\widehat{\delta}=1$
\[
\omega_{n-1}\int\ \log(1/|x|) \,  \vp(x)\, dx= \omega_{n-1}\int p_\mu \, d \delta=\int \big( \widehat{\vp}(\xi) -\widehat{\vp}(0) K_n(|\xi|\big)\, \frac{d\xi}{|\xi|^n},
\]
where the integral is improper in the case $n=1$. Since $(\wh{\vp})^\wedge= (2\pi)^n \vp(- \, \cdot)$, replacement of $\vp$ by $\wh{\vp}$  implies that we have  
\begin{equation}\label{T}
-(2\pi)^{-n}\omega_{n-1}(\log(|\cdot|)^\wedge=T
\end{equation}
 as temperate distribution, where
\[
T\psi :=  \int \big(\psi(\xi) -\psi(0) K_n(|\xi|\big)\, \frac{d\xi}{|\xi|^n}\\
= \lim_{\varepsilon \to 0}\Big(\int_{\varepsilon \le |\xi|} \psi(\xi)\frac{d\xi}{|\xi|^n}-d_n(\varepsilon) \psi(0)\Big)\quad (\psi \in \mathcal{S})
\]
with 
\[
d_n(\varepsilon) :=\omega_{n-1}\int_\varepsilon^\infty K_n(r)\frac{dr}{r}
\]
(cf.\ \cite{Schwartz}, p.\ 44, 258). In this way,  $(\log|\cdot|)^\wedge$ is represented in a different form from the more common  
\[
-(2\pi)^{-n}\omega_{n-1}(\log|\cdot|)^\wedge=S_s+ c_n(s)\delta
\]
 with 
\[
S_s\psi=\int \big(\psi(\xi) -\psi(0) 1_{[0,s]}(|\xi|)\big)\, \frac{d\xi}{|\xi|^n}  \qquad (\psi \in \mathcal{S})
\] 
for $s>0$ and $c_n(s)$  a suitable constant (see e.g.\ \cite{Donoghue}, p.\ 161).   By comparing the two representations it turns out that
\[
c_n(s)= \omega_{n-1}\int_0^\infty( 1_{[0,s]}(r)-K_n(r))\, \frac{dr}{r}\, .
\]
In Section \ref{s3} the integrals are calculated inductively from their values for $n=1,2$.  
\begin{Remark}\label{Rem1}
From the main theorem in \cite{Mattner} it follows that for measures $\mu, \nu \in \mathcal{M}_{+}(\R^n)$  each of the conditions $I(\mu),I(\nu) \in (-\infty, \infty]$ and  $I(\mu),I(\nu) \in [-\infty, \infty)$ is sufficient for the existence of $\int p_\mu\, d\nu$
in  $[-\infty,+\infty]$ and that   
\[
2\int p_\mu\, d\nu \le I(\mu)+I(\nu)
\]
with the integral on the left hand side $>-\infty$  if $I(\mu), I(\nu) \in \R$. In particular, this implies that for  arbitrary $\mu, \nu \in \mathcal{M}_{0}(\R^{n})$ we have
$$\int\int \left\vert\ln\left(\vert x-y\vert\right)\right\vert \, d\vert\mu\vert(y) \, d\vert\nu\vert(x)< \infty\, .
$$
 Interesting enough, $\int p_\mu \, d\nu=-\infty$ may exist for measures $\mu, \nu \in \mathcal{M}_+(\R^n) \setminus \mathcal{M}_{0,+}(\R^n)$, as the example in \cite{Mattner}, p.\ 3340, shows.   
\end{Remark}

In Section \ref{s2} we give the proof of  Theorem \ref{Energy formula several dimensions} (which uses only classical methods) and Section  \ref{s3}  contains some  applications. Applying similar  techniques, in Section  \ref{s4} we present an alternative proof of a known Fourier integral formula for Riesz energy. 

\section{Proof of Theorem \ref{Energy formula several dimensions}}\label{s2}
 
  A main goal is that the formulas for nonnegative measures  hold both in the cases that the energy is finite and is not. We start with three quite simple tools, where the first one is an "improper" version of the Cauchy-Frullani theorem (cf.\ \cite{Ost}).  Recall that $L_{\textup{loc}}(0,\infty)$ is the set of all Lebesgue measurable functions $g:(0,\infty)\rightarrow\C$ such that $g1_{K}$ is integrable for all compact $K\subset(0,\infty)$. Moreover, we write $\|g\|_\infty$ for the essential supremum of $|g|$.

\begin{Lemma} \label{Lemma g 1} Let $g \in L_\textup{loc}(0, \infty)$ and $u \in \C$.  If $t\mapsto g(t)/t$ is improperly integrable at $\infty$ and  $t \mapsto (g(t)-u)/t$ is improperly integrable at $0$, then for arbitrary $a>0$ 
$$
\lim\limits_{\substack{\varepsilon\rightarrow 0\\ N\rightarrow\infty}} \displaystyle\int_{\varepsilon}^{N} \dfrac{g(ra)-g(r)}{r} \, dr= u\ln\left(\dfrac{1}{a}\right)\, .
$$
If, in addition, $g$ is essentially bounded, then
\[
\left\vert\displaystyle\int_{\varepsilon}^{N} \dfrac{g(ra)-g(r)}{r} \, dr\right\vert \leq 2( \|g\|_\infty+ |u|)|\ln(a)|.
\]
\end{Lemma}
\begin{proof} We may assume that $0<a<1$. For fixed $0<\varepsilon<N<+\infty$ we have
\begin{eqnarray*}
\int_{\varepsilon}^{N} \frac{g(ra)-g(r)}{r} \, dr &
= &  \int_{\varepsilon a}^{N a}\frac{g(r)}{r}\, dr - \int_{\varepsilon}^N \frac{g(r)}{r}\, dr
 =  \left(\int_{\varepsilon a}^{\varepsilon}- \int_{N a}^N\right) \frac{g(r)}{r}\, dr\\
& = & u \int_{\varepsilon a}^{\varepsilon}\frac{dr}{r} + \int_{\varepsilon a}^{\varepsilon}\frac{g(r)-u}{r}\, dr - \int_{N a}^N \frac{g(r)}{r}\, dr\; .
\end{eqnarray*}
The first term is $u \ln(1/a)$. According to Cauchy's criterion, the second term  vanishes  as  $\varepsilon $ tends to $0$ and the third as $N$ tends to $\infty$. Moreover, if $g$ is essentially bounded  we obtain that
$$
\left\vert\displaystyle\int_{\varepsilon}^{N} \dfrac{g(ra)-g(r)}{r} \, dr\right\vert \leq 2(\|g\|_\infty + |u|) \left\vert\ln(a)\right\vert.
$$\; .
 \end{proof}
 \begin{Lemma} \label{Lemma g 2} 
 Let $g \in L_\mathrm{loc}(0, \infty)$ be real-valued with $u:=\esssup g<\infty$ and such that $t \mapsto (u-g(t))/t$ is integrable at $0$. 
Then, for all $N>0$ we have
$$
\displaystyle\int_{0}^{N} \dfrac{g(ra)-g(r)}{r} \, dr \begin{cases} \geq 0, & a\in [0,1]\\ \leq 0, & a\in[1,\infty) \end{cases}.
$$
\end{Lemma} 
\begin{proof} The assertion follows from the (essential) nonnegativity of $u-g$ and   
\begin{eqnarray*}
\displaystyle\int_{0}^{N} \dfrac{g(ra)-g(r)}{r} \, dr&=& -\displaystyle\int_{0}^{N} \dfrac{u-g(ra)}{r} \, dr+ \displaystyle\int_{0}^{N} \dfrac{u-g(r)}{r} \, dr\\
&=& -\displaystyle\int_{0}^{aN} \dfrac{u-g(r)}{r} \, dr+ \displaystyle\int_{0}^{N} \dfrac{u-g(r)}{r} \, dr.
\end{eqnarray*}
\end{proof}
\begin{Lemma} \label{Lemma Fatou -infty} Let $g$ satisfy the assumptions of Lemma \ref{Lemma g 1} and Lemma \ref{Lemma g 2} with the same $u$. Then, there are $R,c >0$ such that for each $a>R$ and each $0<\varepsilon<1<N$
$$
\displaystyle\int_{\varepsilon}^{N} \dfrac{g(at)-g(t)}{t} \, dt\leq \displaystyle\int_{0}^{1} \dfrac{u-g(t)}{t} \, dt+c.
$$
\end{Lemma}
\begin{proof} We choose $R>0$ such that
$$
\left\vert\displaystyle\int_{x}^{y} \dfrac{g(t)}{t} \, dt\right\vert \leq 1
$$
whenever $x,y\geq R.$ If $a\in[R,\infty)$ and $N>1>\varepsilon>0,$ then
\begin{eqnarray*}
\displaystyle\int_{\varepsilon}^{N} \dfrac{g(at)-g(t)}{t} \, dt&=& \displaystyle\int_{\varepsilon}^{1} \dfrac{g(at)-g(t)}{t} \, dt+ \displaystyle\int_{1}^{N} \dfrac{g(at)}{t} \, dt - \displaystyle\int_{1}^{N} \dfrac{g(t)}{t} \, dt\\
&=& \displaystyle\int_{\varepsilon}^{1} \dfrac{g(at)-u}{t} \, dt + \displaystyle\int_{\varepsilon}^{1} \dfrac{u-g(t)}{t} \, dt \\
&+& \displaystyle\int_{a}^{aN} \dfrac{g(t)}{t} \, dt- \displaystyle\int_{1}^{N} \dfrac{g(t)}{t} \, dt\\
&\leq& \displaystyle\int_{0}^{1} \dfrac{u-g(t)}{t} \, dt + c
\end{eqnarray*}
where the constant $c>0$ is chosen such that
$\underset{N>1}{\textup{sup}} \, \displaystyle\int_{1}^{N} t^{-1}g(t) \, dt \geq -c+1$.
\end{proof}

\begin{Remark}\label{Remark Kn}
With the aid of the functions  $K_n$ one can express the Fourier transform of the surface measure $\sigma_{n-1}$. 
More precisely, since 
\[
\omega_{n-1}=2\pi^{n/2}/\Gamma(n/2)
\]
for $n\in\N$, we have
(see e.g.\ \cite{Stein2}, p. 154, or \cite{Grafakos}, p. 428)
\begin{equation}\label{Fourier transform sphere}
\widehat{\sigma}_{n-1}(\xi)= \omega_{n-1} K_n(|\xi|)\quad (\xi\in\R^{n}).
\end{equation}
In particular, 
\[
-1\le K_{n}(t)\leq \dfrac{1}{\omega_{n-1}}\displaystyle\int_{\mathbb{S}^{n-1}} 1 \, d\sigma_{n-1}(\xi)=K_{n}(0)=1 \quad (t \in \R) .
\]
Due to the fact that $J_{\alpha}(t)=O(t^{-1/2})$ as $t\rightarrow \infty$ for  $\alpha\geq 0$ (see \cite{Olver}, 10.17.3) one sees that $K_n|_{(0,\infty)}$
satisfies the assumptions of Lemma \ref{Lemma g 1} with $u=1$  and,  
since $K_{n}(t)=1+O(t^{2})$ as $t\rightarrow 0$,  also the assumptions of Lemma \ref{Lemma g 2} are fulfilled.
Finally,  because $K_{n}'(0)=0$ and $K_{n}''(0)<0$ we see that $K_{n}$ is decreasing in some interval $[0,\delta].$  Writing $ \Delta_a(t):=K_n(at)-K(t)$, for $a\in[0,1)$ and $N>1>\delta>\varepsilon>0$ we thus obtain by Lemma \ref{Lemma g 2} 
\begin{eqnarray*}
\displaystyle\int_{\varepsilon}^{N} \dfrac{K_{n}(at)-K_{n}(t)}{t} \, dt &=& \displaystyle\int_{0}^{N} \dfrac{\Delta_a(t)}{t} \, dt- \displaystyle\int_{0}^{\varepsilon} \dfrac{\Delta_a(t)}{t} \, dt
\geq - \displaystyle\int_{0}^{\varepsilon} \dfrac{\Delta_a(t)}{t} \, dt\\
&\geq & - \displaystyle\int_{0}^{\delta} \dfrac{\Delta_a(t)}{t} \, dt
\geq - \displaystyle\int_{a\delta}^{a} \dfrac{\Delta_a(t)}{t} \, dt \\
&\geq & - \displaystyle\int_{0}^{1} \dfrac{1-K_{n}(t)}{t} \, dt\, .
\end{eqnarray*}
\end{Remark}
Now, we are in a position to give the
\begin{proof}[Proof of Theorem \ref{Energy formula several dimensions}:]
Fix $0<\varepsilon<N<+\infty.$ Then, by Fubini's theorem and  \eqref{Lemma sphere} 
\begin{align*}
&\displaystyle\int_{\varepsilon\leq \vert \xi\vert\leq N} \left(\widehat{\mu}(\xi)\cdot\overline{\widehat{\nu}(\xi)}-\mu(\R^{n})\bar{\nu}(\R^{n})K_{n}(\vert \xi\vert)\right) \, \dfrac{d\xi}{\vert \xi\vert^{n}}\\
&= \displaystyle\int_{\varepsilon\leq \vert \xi\vert\leq N} \left[\left(\displaystyle\int e^{i \xi x} \, d\mu(x)\right)\cdot\left(\displaystyle\int e^{-i\xi  y} \, d\bar{\nu}(y)\right)-\mu(\R^{n})\bar{\nu}(\R^{n})K_{n}(\vert \xi\vert)\right] \, \dfrac{d\xi}{\vert \xi\vert^{n}}\\
&= \displaystyle\int\displaystyle\int_{\varepsilon\leq \vert \xi\vert\leq N} \frac{e^{i\xi(x-y)}-K_{n}(\vert \xi\vert)}{\vert \xi\vert^{n}} \, d\xi \, d(\mu\otimes\bar{\nu})(x,y)\\
&= \displaystyle\int\displaystyle\int_{\mathbb{S}^{n-1}} \displaystyle\int_{\varepsilon}^{N} \frac{e^{ir\zeta(x-y)}-K_{n}(r)}{r} \, dr \, d\sigma_{n-1}(\zeta) \, d(\mu\otimes\bar{\nu})(x,y).
\end{align*}
Another application of Fubini's theorem and  \eqref{Fourier transform sphere} give us
\begin{align*}
&\displaystyle\int\displaystyle\int_{\mathbb{S}^{n-1}} \displaystyle\int_{\varepsilon}^{N} \frac{e^{ir\zeta(x-y)}-K_{n}(r)}{r} \, dr \, d\sigma_{n-1}(\zeta) \, d(\mu\otimes\bar{\nu})(x,y)\\
&=\displaystyle\int \displaystyle\int_{\varepsilon}^{N} \left[ \left(\, \displaystyle\int_{\mathbb{S}^{n-1}} e^{ir(x-y)} \, d\sigma_{n-1}(\zeta)\right) -\omega_{n-1}K_{n}(r)\right] \, \dfrac{dr}{r} \, d(\mu\otimes\bar{\nu})(x,y)\\
&=\omega_{n-1}\displaystyle\int \displaystyle\int_{\varepsilon}^{N} \dfrac{K_{n}(r\vert x-y\vert)-K_{n}(r)}{r} \, dr \, d(\mu\otimes\bar{\nu})(x,y).
\end{align*}
First, suppose that
$$
\displaystyle\int\displaystyle\int\vert\ln(\vert x-y\vert)\vert \, d\vert\mu\vert(x) \, d\vert\nu\vert(y) <+\infty.
$$
Then we necessarily have $(\vert\mu\vert\otimes\vert\nu\vert)(\{(x,y)\in\R^{n}\times\R^{n}: x=y\})=0.$ Therefore, if we take $g=K_{n}$ and $u=g(0)$ in Lemma \ref{Lemma g 1} (see Remark \ref{Remark Kn}) the dominated convergence theorem yields that
\begin{align*}
&\lim\limits_{\substack{\varepsilon\rightarrow 0\\ N\rightarrow\infty}} \displaystyle\int \displaystyle\int_{\varepsilon}^{N} \dfrac{K_{n}(r\vert x-y\vert)-K_{n}(r)}{r} \, dr \, d(\mu\otimes\bar{\nu})(x,y)\\
&= \displaystyle\int \ln\left(\dfrac{1}{\vert x-y\vert}\right) \, d(\mu\otimes\bar{\nu})(x,y) = \int p_\mu \, d\bar{\nu}.
\end{align*}
Now, let $\mu,\nu\in\mathcal{M}_{+}(\R^{n})$ be such that
$$
\displaystyle\int\ln\left(\dfrac{1}{\vert x-y\vert}\right) \, d(\mu\otimes\nu)(x,y) =+\infty.
$$
In this case, 
$$
\displaystyle\int_{\{\vert x-y\vert<1\}} \ln\left(\dfrac{1}{\vert x-y\vert}\right) \, d(\mu\otimes\nu)(x,y) =+\infty
$$
since 
$$
\displaystyle\int_{\{\vert x-y\vert\geq 1\}} \ln\left(\dfrac{1}{\vert x-y\vert}\right) \, d(\mu\otimes\nu)(x,y)\in (-\infty,0].
$$
The same argument implies together with Lemma \ref{Lemma g 1}, Remark \ref{Remark Kn} and the dominated convergence theorem that
$$
\lim\limits_{\substack{\varepsilon\rightarrow 0\\ N\rightarrow\infty}} \ \displaystyle\int_{\{\vert x-y\vert\geq 1\}}\left(\displaystyle\int_{\varepsilon}^{N} \dfrac{K_{n}(r\vert x-y\vert)-K_{n}(r)}{r} \, dr\right) \, d(\mu\otimes\nu)(x,y)
$$
exists and is finite. Therefore, we only have to show that
$$
\liminf\limits_{\substack{\varepsilon\rightarrow 0\\ N\rightarrow\infty}} \displaystyle\int_{\{\vert x-y\vert< 1\}}\left(\displaystyle\int_{\varepsilon}^{N} \dfrac{K_{n}(r\vert x-y\vert)-K_{n}(r)}{r} \, dr\right) \, d(\mu\otimes\nu)(x,y)= +\infty.
$$
By Remark \ref{Remark Kn}, we know that
$$
\displaystyle\int_{\varepsilon}^{N} \dfrac{K_{n}(r\vert x-y\vert)-K_{n}(r)}{r} \, dr
\geq - \displaystyle\int_{0}^{1} \dfrac{1-K_{n}(r)}{r} \, dr
$$
for $\vert x-y\vert<1$ and sufficiently large $N$ and small $\varepsilon.$ Therefore, we can apply Fatou's lemma (notice that $\mu\otimes\nu$ is finite) and get with $K_n(0)=1$
\begin{align*}
&\liminf\limits_{\substack{\varepsilon\rightarrow 0\\ N\rightarrow\infty}} \displaystyle\int_{\{\vert x-y\vert< 1\}}\left(\displaystyle\int_{\varepsilon}^{N} \dfrac{K_{n}(r\vert x-y\vert)-K_{n}(r)}{r} \, dr\right) \, d(\mu\otimes\nu)(x,y)\\
&\geq \displaystyle\int_{\{x=y\}}\liminf\limits_{\substack{\varepsilon\rightarrow 0\\ N\rightarrow\infty}} \left(\displaystyle\int_{\varepsilon}^{N} \dfrac{1-K_{n}(r)}{r} \, dr\right) \, d(\mu\otimes\nu)(x,y)\\ &+ \displaystyle\int_{\{0<\vert x-y\vert< 1\}}\ln\left(\dfrac{1}{\vert x-y\vert}\right) \, d(\mu\otimes\nu)(x,y)
\end{align*}
If $(\mu\otimes\nu)(\{(x,y)\in\R^{n}\times\R^{n}: x=y\})>0,$ then the first integral is equal to $+\infty$ since
$$
\liminf\limits_{\substack{\varepsilon\rightarrow 0\\ N\rightarrow\infty}} \displaystyle\int_{\varepsilon}^{N} \dfrac{1-K_{n}(t)}{t} \, dt= +\infty.
$$
If $(\mu\otimes\nu)(\{(x,y)\in\R^{n}\times\R^{n}: x=y\})=0,$  the second integral is equal to $+\infty.$

Suppose finally that $\mu,\nu\in\mathcal{M}_{+}(\R^{n})$ and 
$$
\displaystyle\int\ln\left(\dfrac{1}{\vert x-y\vert}\right) \, d(\mu\otimes\nu)(x,y)=-\infty.
$$
Then, 
$$
\displaystyle\int_{\{\vert x-y\vert\geq T\}} \ln\left(\dfrac{1}{\vert x-y\vert}\right) \, d(\mu\otimes\nu)(x,y)=-\infty
$$
for all $T\geq 1$ since
$$
\displaystyle\int_{\{\vert x-y\vert\leq T\}} \ln\left(\dfrac{1}{\vert x-y\vert}\right) \, d(\mu\otimes\nu)(x,y)\in\R.
$$
The same argument implies together with Lemma \ref{Lemma g 1} and the dominated convergence theorem that
$$
\lim\limits_{\substack{\varepsilon\rightarrow 0\\ N\rightarrow\infty}} \ \displaystyle\int_{\{\vert x-y\vert\leq T\}}\left(\displaystyle\int_{\varepsilon}^{N} \dfrac{K_{n}(r\vert x-y\vert)-K_{n}(r)}{r} \, dr\right) \, d(\mu\otimes\nu)(x,y)
$$
exists and is finite for all $T\geq 1.$  
If we pick $g=K_{n}$ in Lemma \ref{Lemma Fatou -infty}, then there is some $R>1$ such that for all $\vert x-y\vert\geq R$ and $0<\varepsilon<1<N$ 
$$
\displaystyle\int_{\varepsilon}^{N} \dfrac{K_{n}(r\vert x-y\vert)-K_{n}(r)}{r} \, dr\leq \displaystyle\int_{0}^{1} \dfrac{1-K_{n}(r)}{r} \, dr + c
$$
for some constant $c\in\R$ independent of $x,y,\varepsilon$ and $N.$
Therefore, we can apply Fatou's lemma since $\mu\otimes\nu$ is finite and get
\begin{align*}
&\limsup\limits_{\substack{\varepsilon\rightarrow 0\\ N\rightarrow\infty}} \displaystyle\int_{\{\vert x-y\vert\geq R\}}\left(\displaystyle\int_{\varepsilon}^{N} \dfrac{K_{n}(r\vert x-y\vert)-K_{n}(r)}{r} \, dr\right) \, d(\mu\otimes\nu)(x,y)\\
&\leq \displaystyle\int_{\{\vert x-y\vert\geq R\}}\ln\left(\dfrac{1}{\vert x-y\vert}\right) \, d(\mu\otimes\nu)(x,y)=-\infty.
\end{align*}
 \end{proof}

\section{Consequences of Theorem \ref{Energy formula several dimensions}} \label{s3}

Let us start with two illustrating examples, which also show that  some interesting  integrals involving $J_0$ emerge form Theorem \ref{Energy formula several dimensions}: 

According  to \eqref{Fourier transform sphere},  Theorem \ref{Energy formula several dimensions} and \eqref{Lemma sphere} imply  that
\[
I(\sigma_{n-1})=\omega_{n-1}^2\int_0^\infty  (K_n^2(r) -K_n(r))\frac{dr}{r}\qquad (n \in \N) .
\]
Since $K_2=J_0$ and since, as already mentioned in the introduction, $I(\sigma_1)= 0$, we obtain  
\begin{equation}\label{Im0}
0=\int_0^\infty  (J_0^2(r) -J_0(r))\frac{dr}{r}\, .
\end{equation}
Applying Lemma \ref{Lemma g 1} with $g=J_0^2$ we get, more generally, 
\begin{equation}\label{Im}
\int_0^\infty (J_0^2(ar) -J_0(r))\frac{dr}{r}=\ln \left(\frac{1}{a}\right)
\end{equation}
for $a>0$. 
The standard one-dimensional example in which $\wh{\mu}$ is known is the arcsine distribution 
\[
d\mu(t)=\frac{dt}{\pi\sqrt{1-t^2}} 
\]
on $[-1,1]$. Here is (see \cite{Erdelyi}, p. 11) $J_0=\wh{\mu}$, but then with the Dirac measure $\delta_0$ also
\[
(\mu\otimes \delta_0)^{\wedge}(\xi_1, \xi_2)=\wh{\mu}(\xi_1)=J_0(\xi_1) \qquad (\xi_1, \xi_2\in \R).
\]  
Since $I(\mu)=I(\mu \otimes \delta_0)$, by choosing $n=2$ in Theorem \ref{Energy formula several dimensions}   and using \eqref{Im} we have
\[
I(\mu)= \frac{1}{2\pi}\int_{-\pi}^\pi \int_0^\infty (J_0^2(r \cos\theta)-J_0(r)) \frac{dr}{r}\, d \theta
= \frac{2}{\pi} \int_0^{\pi/2} \ln\left(\frac{1}{\cos \theta}\right) d\theta =  \ln 2\, .
\]
This is of course folklore but usually proved in a quite different way by  exploiting some amount of  potential theory in the plane. More precisely, one shows that  $\mu$ is the equilibrium measure of $[-1,1]$, that is, 
$\mu$ minimizes logarithmic energy among all Borel probability measures on $[-1,1]$, and that the minimum is $\ln 2$ or, in other words, that the logarithmic capacity of $[-1,1]$ is $1/2=e^{-\ln 2}$.

Also, taking $n=1$ in Theorem \ref{Energy formula several dimensions} and using \eqref{Im0} leads to  
\[
\ln 2= \int_0^\infty (J_0^2(r)-\cos r)\,\frac{dr}{r}= \int_0^\infty (J_0(r)-\cos r)\,\frac{dr}{r}
\]
(see e.g.\  \cite{BGKM}, Example 5.5,  \cite{Humbert}, p.\ 278). According to Lemma \ref{Lemma g 1}, we end at 
\begin{equation} \label{Humbert2}
\int_0^\infty (J_0(ar)-\cos r)\frac{dr}{r}=\ln \left(\frac{2}{a}\right) 
\end{equation}
for arbitrary $a>0$. In particular, this implies that  for $n=1$ the function $K_1=\cos$ in Theorem \ref{Energy formula several dimensions} may be replaced by $\xi \mapsto J_0(2|\xi|)$.

\begin{Remark}  \label{cn(s)} (cf. \cite{Schwartz}, p.\ 258, and  \cite{Vladimirov}, p.\ 118, for the case $n=2$)
Let $s>0$. It is known that
\[
\displaystyle\int_{0}^{\infty} (1_{[0,s]}(r)-\cos(r)) \, \frac{dr}{r} = \gamma+\ln(s)\, ,
\]
where $\gamma$ is the Euler-Mascheroni constant (see e.g.\ \cite{Olver}, 6.2.13). In combination with (\ref{Humbert2}) this gives 
\[
\int_{0}^{\infty} (1_{[0,s]}(r)-J_{0}(r))\, \frac{dr}{r}= \gamma+\ln(s)-\ln(2)
\]
More generally, we have
\[
\displaystyle\int_{0}^{\infty} (1_{[0,s]}(r)-K_{n}(r))\, \frac{dr}{r}= \begin{cases} \gamma+\ln(s/2)-\sum\limits_{k=1}^{n/2-1} \frac{1}{2k}, & n \mbox{ even }\\  \gamma+\ln(s)-\sum\limits_{k=1}^{(n-1)/2} \frac{1}{2k-1}, & n \mbox{ odd }\end{cases}.
\]
Indeed:  Since $K_{1}=\cos$ and $K_{2}=J_0$ the formula holds for $n=1$ and $n=2$. Now, suppose that the formula is true for  $n\in\N.$ The recurrence formula for Bessel functions (see \cite{Olver}, 10.6.1) yields that
\[
J_{n/2}(t)=\frac{t}{n}(J_{n/2-1}(t)+J_{n/2+1}(t)) \quad (t>0).
\]
Since 
\[
\int_0^\infty r^{-n/2} J_{n/2+1}(r)\, dr = \frac{1}{2^{n/2}\Gamma(n/2+1)}
\]
 (see e.g. \cite{Olver}, 10.22.43) we obtain that
\begin{eqnarray*}
\displaystyle\int_{0}^{\infty} (1_{[0,s]}(r)-K_{n+2}(r))\, \frac{dr}{r}
&=& \displaystyle\int_{0}^{\infty} \left(1_{[0,s]}(r)-\frac{2^{n/2}\Gamma\left(n/2+1\right)J_{n/2}(r)}{r^{n/2}}\right)\, \frac{dr}{r}\\
&=& \displaystyle\int_{0}^{\infty} (1_{[0,s]}(r)-K_{n}(r))\, \frac{dr}{r} - \frac{1}{n} 
\end{eqnarray*}
and hence the formula holds for $n+2$.

Since 
\[
-\int_\eps^\infty K_n(r)\frac{dr}{r}= \ln(\eps) +\int_{\eps}^{\infty} (1_{[0,1]}(r)-K_{n}(r))\, \frac{dr}{r},
\]
according to \eqref{T} we recover the representation of $(\log |\cdot|)^\wedge$ from \cite{Schwartz}, p.\ 258. 
\end{Remark}
A next  consequence of Theorem \ref{Energy formula several dimensions} is a characterization for the finiteness of the logarithmic energy for nonnegative measures.
\begin{Corollary} \label{Energy Fourier several dimensions} Let $\mu\in\mathcal{M}_{0, +}(\R^{n})$ be such that
\[
\lim_{\varepsilon\to 0} \int_{\varepsilon\le |\xi|\le 1} (|\wh{\mu}(\xi)|^2-|\mu(\R^n)|^2)\frac{d\xi}{|\xi|^n}
\]
exists in $\R$.
Then $I(\mu)$ is finite if and only if $\vert\widehat{\mu}|^2$ is locally integrable at $\infty$ with respect to $d\xi/|\xi|^n$.
\end{Corollary}
\begin{proof} By Theorem \ref{Energy formula several dimensions} we know that
$$
\omega_{n-1}I(\mu)=\lim\limits_{\substack{\varepsilon\rightarrow 0\\ N\rightarrow\infty}} \ \displaystyle\int_{\varepsilon\leq \vert \xi\vert\leq N} \left(\vert\widehat{\mu}(\xi)\vert^{2}-\vert\mu(\R^{n})\vert^{2}K_{n}(\vert \xi\vert)\right) \, \dfrac{d\xi}{\vert \xi\vert^{n}}.
$$
According to the assumption and the smoothness of $K_n$, the limit in the equation above exists if and only if
$$
\lim\limits_{N\rightarrow\infty} \displaystyle\int_{1\leq \vert \xi\vert\leq N}  \left(\vert\widehat{\mu}(\xi)\vert^{2}-\vert\mu(\R^{n})\vert^{2}K_{n}(\vert \xi\vert)\right) \, \dfrac{d\xi}{\vert \xi\vert^{n}}
$$
exists.
Recalling that
$$
\lim\limits_{N\rightarrow\infty} \displaystyle\int_{1\leq \vert \xi\vert\leq N} \dfrac{K_{n}(\vert \xi\vert)}{\vert \xi\vert^{n}} \, d\xi=\omega_{n-1} \lim\limits_{N\rightarrow\infty} \displaystyle\int_{1}^{N} \dfrac{K_{n}(r)}{r} \, dr 
$$
is finite,  we conclude that $I(\mu)\in\R$ if and only if
$$
\lim\limits_{N\rightarrow\infty} \displaystyle\int_{1\leq \vert \xi\vert\leq N} \dfrac{\vert\widehat{\mu}(\xi)\vert^{2}}{\vert \xi\vert^{n}} \, d\xi= \displaystyle\int_{\vert \xi\vert\geq 1} \dfrac{\vert\widehat{\mu}(\xi)\vert^{2}}{\vert \xi\vert^{n}} \, d\xi<\infty. 
$$
\end{proof}

As already mentioned in the introduction, the limit $\varepsilon \to 0$ in the preceding corollary  exists for arbitrary Borel measures having compact support. As a consequence, the equivalence statement of the corollary  also holds for compactly supported complex measures.

\begin{Remark}
We consider the linear space 
$\mathcal{M}_{00}(\R^n)$ of all $ \mu \in \mathcal{M}_{0}(\R^n)$ with vanishing total mass, that is 
$\mu(\R^n)=0$.  For $\mu \in \mathcal{M}_{00}(\R^n)$ and  arbitrary $\nu \in  \mathcal{M}_{0}(\R^n)$,  Theorem \ref{Energy formula several dimensions} says that
\[
\omega_{n-1} \int p_\mu\, d\bar{\nu}= \lim\limits_{\substack{\varepsilon\rightarrow 0\\ N\rightarrow\infty}} \ \displaystyle\int\limits_{\varepsilon\leq \vert \xi\vert\leq N}(\widehat{\mu}\, \overline{\widehat{\nu}})(\xi) \, \dfrac{d\xi}{\vert \xi\vert^{n}}.
\]
In particular, for $\nu=\mu$ we get  
\[
\omega_{n-1} I(\mu)  = \int  |\widehat{\mu}(\xi) |^2\, \dfrac{d\xi}{\vert \xi\vert^{n}} \in [0, \infty) \, .
\]
This may be seen as a continuous version of \eqref{discrete}. It implies in particular that the energy integral is positive definite on the space $\mathcal{M}_{00} (\R^n).$ This fact is well known in the case of compactly supported, signed (real) measures in $\mathcal{M}_{00}(\R^n)$ (see  \cite{Garnett3}, Chapter III, Theorem 3.1, \cite{ST}, Chapter I, Lemma 1.8, cf.\ also \cite{Fukushima}, \cite{Deny}, 5, III.5). Moreover (see \cite{Mattner}), if $\mu$ is a not necessarily finite signed  measure with $\mu(\R^n)=0$ and 
$$
\displaystyle\int\displaystyle\int \ln\left(\dfrac{1}{\vert x-y\vert}\right) \, d\mu(x) \, d\mu(y)
$$
exists in $[-\infty,+\infty],$ then this integral is either $\geq 0$ or equal to $+\infty.$
\end{Remark}
Since, according to  continuity,  $|\wh{\mu}|^2$ is not locally integrable with respect to $d\xi/|\xi|^n$ at the origin if $\mu(\R^n)\not=0$, we have 

\begin{Corollary} Let $\mu\in\mathcal{M}_{0}(\R^{n})$. Then  $\wh{\mu} \in L^2(d\xi/|\xi|^n)$ if and only if  $\mu \in\mathcal{M}_{00}(\R^{n})$  and in this case 
$\omega_{n-1} I(\mu)=\|\wh{\mu}\|^2_{L^2(d\xi/|\xi|^n)}$.
\end{Corollary} 

\section{Riesz Energy}\label{s4}

Logarithmic energy may be viewed as the limit case $\alpha \to 0$ of Riesz energies. For positive $\alpha$ and $x \in \R^n$ we write
$$
p_{\mu,\alpha}(x):=\displaystyle\int \dfrac{1}{\vert x-y\vert^{\alpha}} \, d\mu(y)\, \in \C \cup \{+\infty\} 
$$
if $y \mapsto |x-y|^{-\alpha}$ is integrable with respect to $|\mu|$ or $\mu \in \mathcal{M}_+(\R^n)$.
Moreover, we write  $\mathcal{M}_{\alpha}(\R^{n})$ for the set of all complex Borel measures with 
$$
\displaystyle\int\displaystyle\int \dfrac{1}{\vert x-y\vert^{\alpha}} \, d\vert\mu\vert(x) \, d\vert\mu\vert(y)<+\infty\; .
$$
For $\mu\in\mathcal{M}_{\alpha}(\R^{n})\cup \mathcal{M}_+(\R)$ the Riesz potential $p_{\mu, \alpha}$ is defined $\mu$ almost everywhere  and  
$$
I_{\alpha}(\mu):= \displaystyle\int p_{\mu,\alpha}(x) \, d\bar{\mu}(x)= \displaystyle\int\displaystyle\int \dfrac{1}{\vert x-y\vert^{\alpha}} \, d\mu(y) \, d\bar{\mu}(x) 
$$
is called the Riesz energy of order $\alpha$ of $\mu$.
It is a well-known fact (see \cite{Mattila}, p. 162,  \cite{Landkof}, Chapter VI) that if $\mu\in\mathcal{M}_+(\R^n)$ has compact support or $\mu\in\mathcal{M}_{\alpha}(\R^n)$ and $\alpha<n$, then
\begin{equation}\label{mattila}
\gamma_{n,\alpha}I_{\alpha}(\mu)= \displaystyle\int \dfrac{\vert\widehat{\mu}(\xi)\vert^{2}}{\vert \xi\vert^{n-\alpha}} \, d\xi\, , 
\end{equation}
where  
\[
\gamma_{n,\alpha}=\omega_{n-1}\dfrac{2^{\alpha-1}\Gamma\left(\frac{n}{2}\right)\Gamma\left(\frac{\alpha}{2}\right)}{\Gamma\left(\frac{n-\alpha}{2}\right)}= \dfrac{2^{\alpha}\pi^{\frac{n}{2}}\Gamma\left(\frac{\alpha}{2}\right)}{\Gamma\left(\frac{n-\alpha}{2}\right)}\, .
\]
 Here the situation is more convenient compared to the logarithmic case due to the fact that  $\xi \mapsto|\xi|^\alpha |\widehat{\mu}|^2(\xi)$ is locally integrable with respect to $d\xi/|\xi|^n$ at the origin. 
Arbitrary complex measures are considered  in \cite{HR1}, where Riesz energy is defined by the formula \eqref{mattila}. Finally, according to the pre-Hilbert space structure of the space of signed measures in $\mathcal{M}_\alpha(\R^n)$ endowed with the inner product 
$(\mu, \nu) \mapsto \int p_{\mu, \alpha} \, d \nu$ (see e.g.\ \cite{Landkof}, p.\ 82), the mixed integral $\int p_{\mu, \alpha}\, d\bar{\nu}$ exists in $\C$ whenever $\mu, \nu \in \mathcal{M}_\alpha(\R^n)$ (and of course in $[0, \infty]$ for nonnegative $\mu, \nu$). 
 
Applying similar techniques as in the proof of Theorem \ref{Energy formula several dimensions} one can show the following result, in which $\alpha$ is restricted  to be less than  $(n+1)/2$ due to the fact that  $K_n$ has to be (at least improperly) integrable at $\infty$ with respect to $t^{\alpha-1}dt$.  
 
\begin{theorem} \label{Riesz}
Let $0<\alpha<(n+1)/2$.  If $\mu, \nu \in \mathcal{M}_+(\R^n)$ or if $\mu, \nu \in \mathcal{M}(\R^{n})$ with 
$$
\displaystyle\int\displaystyle\int \dfrac{1}{\vert x-y\vert^{\alpha}} \, d\vert\mu\vert(x) \, d\vert\nu\vert(y)<+\infty
$$
then 
$$
\gamma_{n,\alpha}\int p_{\mu, \alpha} \,  d\bar{\nu}=\lim\limits_{N\rightarrow\infty} \displaystyle\int_{\vert \xi\vert\leq N} \vert \xi\vert^{\alpha} \, (\widehat{\mu}\, \overline{\widehat{\nu}})(\xi) \, \frac{d\xi}{|\xi|^n}\, .
$$
In particular, for $\mu \in \mathcal{M}_+(\R^n) \cup \mathcal{M}_\alpha(\R^n)$ we have
$$
\gamma_{n,\alpha}I_{\alpha}(\mu)=\displaystyle\int \vert \xi\vert^{\alpha} \vert\widehat{\mu}(\xi)\vert^{2} \, \frac{d\xi}{|\xi|^n}.
$$
\end{theorem}
\begin{proof} Let $N>0.$ Then, by Fubini's theorem,  \eqref{Lemma sphere} and \eqref{Fourier transform sphere}
\begin{eqnarray*}
\displaystyle\int_{\vert \xi\vert \leq N} \vert \xi\vert^{\alpha} (\widehat{\mu}\,\overline{\widehat{\nu}})(\xi) \, \frac{d\xi}{|\xi|^n}
&=&\displaystyle\int\displaystyle\int \left(\, \displaystyle\int_{\vert \xi\vert\leq N} e^{it(x-y)} \vert t\vert^{\alpha-1} \, dt\right) \, d\mu(x) \, \bar{\nu}(y)\\
&=& \displaystyle\int \left(\, \displaystyle\int \displaystyle\int_{0}^{N} e^{ir\zeta(x-y)} r^{\alpha-1} \, dr \, d\sigma_{n-1}(\zeta)\right)d(\mu\otimes\bar{\nu})(x,y)\\
&=& \omega_{n-1} \displaystyle\int\left(\displaystyle\int_{0}^{N} K_{n}(r\vert x-y\vert)r^{\alpha-1} \, dr\right) \, d(\mu\otimes\bar{\nu})(x,y).
\end{eqnarray*}
First, suppose that 
$$
\displaystyle\int\displaystyle\int\dfrac{1}{\vert x-y\vert^{\alpha}} \, d\vert\mu\vert(x) \, d\vert\nu\vert(y)<+\infty.
$$
Then, we have $(\vert\mu\vert\otimes\vert\nu\vert)(\{(x,y)\in\R^{n}\times\R^n: x=y\})=0.$ If $(x,y)\in\R^{n}$ with $x\neq y,$ then by definition of $K_{n}$ and \cite{Olver}, 10.22.43
\begin{eqnarray*}
\lim\limits_{N\rightarrow\infty} \displaystyle\int_{0}^{N} K_{n}(t\vert x-y\vert) t^{\alpha-1} \, dt &=& \dfrac{2^{\frac{n}{2}-1}\Gamma\left(\frac{n}{2}\right)}{\vert x-y\vert^{\alpha}} \displaystyle\int_0^\infty J_{\frac{n}{2}-1}(t)t^{\alpha-\frac{n}{2}} \, dt\\
&=&\frac{2^{\alpha-1}\Gamma\left(\frac{n}{2}\right)\Gamma\left(\frac{\alpha}{2}\right)}{\Gamma\left(\frac{n-\alpha}{2}\right)\vert x-y\vert^{\alpha}}
\end{eqnarray*}
as well as
$$
\left\vert\displaystyle\int_{0}^{N} K_{n}(\vert x-y\vert t) t^{\alpha-1} \, dt\right\vert\leq \dfrac{1}{\vert x-y\vert^{\alpha}}\cdot \underset{M>0}{\textup{sup}} \left\vert\displaystyle\int_{0}^{M} K_{n}(t) t^{\alpha-1} \, dt\right\vert.
$$
Therefore, the dominated convergence theorem implies that
\begin{align*}
&\lim\limits_{N\rightarrow\infty} \displaystyle\int_{\vert \xi\vert\leq N}  \vert \xi\vert^{\alpha} (\widehat{\mu}\,\overline{\widehat{\nu}})(\xi) \, \frac{d\xi}{|\xi|^n}\\
&= \omega_{n-1} \lim\limits_{N\rightarrow\infty} \displaystyle\int\left(\displaystyle\int_{0}^{N} K_{n}(r\vert x-y\vert)r^{\alpha-1} \, dr\right) \, d(\mu\otimes\bar{\nu})(x,y)\\
&=\gamma_{n,\alpha} \displaystyle\int\displaystyle\int \dfrac{1}{\vert x-y\vert^{\alpha}} \, d\mu(x) \, d\bar{\nu}(y)
\end{align*}
Now, let $\mu$ and $\nu$ be positive measures such that 
$$
\displaystyle\int\dfrac{1}{\vert x-y\vert^{\alpha}} \, d(\mu\otimes\nu)(x,y)=+\infty.
$$
Then, we necessarily have
$$
\displaystyle\int_{\{\vert x-y\vert<1\}} \dfrac{1}{\vert x-y\vert^{\alpha}} \, d(\mu\otimes\nu)(x,y)=+\infty
$$
since 
$$
\displaystyle\int_{\{\vert x-y\vert\geq 1\}} \dfrac{1}{\vert x-y\vert^{\alpha}} \, d(\mu\otimes\nu)(x,y)\in [0,\infty).
$$
The same argument implies together with the dominated convergence that
$$
\lim\limits_{N\rightarrow\infty} \displaystyle\int_{\{\vert x-y\vert\geq 1\}}\left(\displaystyle\int_{0}^{N} K_{n}(r\vert x-y\vert)r^{\alpha-1} \, dr\right) \, d(\mu\otimes\nu)(x,y)
$$
exists and is finite. Therefore, we only have to show that
$$
\liminf\limits_{N\rightarrow\infty} \displaystyle\int_{\{\vert x-y\vert< 1\}}\left(\displaystyle\int_{0}^{N} K_{n}(r\vert x-y\vert)r^{\alpha-1} \, dr\right) \, d(\mu\otimes\nu)(x,y)= +\infty.
$$
But since
\begin{eqnarray*}
\displaystyle\int_{0}^{N} K_{n}(r\vert x-y\vert) r^{\alpha-1} \, dr&=& \dfrac{1}{\vert x-y\vert^{\alpha}} \displaystyle\int_{0}^{\vert x-y\vert N} K_{n}(r) r^{\alpha-1} \, dr\\
&\geq& \underset{M>0}{\textup{inf}} \, \displaystyle\int_{0}^{M} K_{n}(r)r^{\alpha-1} \, dr
\end{eqnarray*}
whenever $\vert x-y\vert<1,$ Fatou's lemma (notice that $\mu\otimes\nu$ is finite) gives us
\begin{align*}
&\liminf\limits_{N\rightarrow\infty} \displaystyle\int_{\{\vert x-y\vert< 1\}}\left(\displaystyle\int_{0}^{N} K_{n}(r\vert x-y\vert)r^{\alpha-1} \, dr\right) \, d(\mu\otimes\nu)(x,y)\\
&\geq \displaystyle\int_{\{x=y\}}\liminf\limits_{N\rightarrow\infty} \left(\displaystyle\int_{0}^{N} r^{\alpha-1} \, dr\right) \, d(\mu\otimes\nu)(x,y)\\
&+ \gamma_{n,\alpha} \displaystyle\int_{\{0<\vert x-y\vert< 1\}}\dfrac{1}{\vert x-y\vert^{\alpha}} \, d(\mu\otimes\nu)(x,y).
\end{align*}
If $(\mu\otimes\nu)(\{(x,y)\in\R^{n}\times\R^n: x=y\})>0,$ then the first integral is equal to $+\infty$ since
\begin{eqnarray*}
\liminf\limits_{N\rightarrow\infty} \displaystyle\int_{0}^{N} r^{\alpha-1} \, dr= +\infty.
\end{eqnarray*}
If $(\mu\otimes\nu)(\{(x,y)\in\R^{n}\times\R^n: x=y\})=0,$ then the second integral equals $+\infty$.
\end{proof}

\vspace{1em}

\textbf{Acknowledgment} 
The authors are grateful to Lutz Mattner, Norbert Ortner, Peter Wagner and Jochen Wengenroth for several helpful suggestions and comments. Special thanks to Norbert Ortner for his letter containing a number of valuable informations around the subject.

\end{document}